\documentclass[12pt]{amsart}
\usepackage{amsmath}
\usepackage{amssymb}
\usepackage{amsfonts}
\usepackage{tikz-cd}

\newtheorem{thm}{Theorem}
\newtheorem{prop}{Proposition}
\newtheorem{lem}{Lemma}
\newtheorem{cor}{Corollary}

\newtheorem{defn}{Definition}
\newtheorem{ex}{Example}
\newtheorem{ack}{Acknowledgement}

\begin{document}
\title{Valdivia's  lifting  theorem for non-metrizable spaces. Preprint.}
\author{Thomas E. Gilsdorf}    
\maketitle
\begin{center} Department of Mathematics \\
 Central Michigan University   \\
 Mt. Pleasant, MI  48859  \/ USA  \\
 gilsd1te@cmich.edu   \\

 \end{center}

\bigskip 

\begin{center}  \today  \end{center} 

\bigskip 

\noindent  {\it Abstract.}  \/   Valdivia's lifting  theorem of (pre) compact sets and convergent (respectively, Cauchy) sequences from a quasi-(LB) space to a metrizable, strictly barrelled space  is extended  to a strictly larger collection of range spaces. Specifically,  we  assume  that the range space has a  sequential web structure and do not require that it be metrizable, nor  strictly barrelled, and the range space need not  even be  barrelled.  Distinguishing examples are provided that include natural constructions of range spaces connected with applications, such as   $\mathcal{D}_{\Gamma}^{\prime}$, the space  of distributions having their wavefront sets in a specified closed cone $\Gamma$.   The same and other examples could also serve as domain spaces for Valdivia's closed graph theorem, revealing  a much wider collection of domain spaces that can be used in that result.     \\
\bigskip

\noindent {\bf 2020 Mathematics Subject Classification:} Primary 46A30; Secondary 54C05, 54C10, 46A08, 47A05, 46A03.  \\

\noindent {\bf Keywords:}  Quasi-(LB) space, strictly barrelled, web, lifting theorem, closed graph theorem.   \\

\section{Introduction.}

 Nearly a decade after DeWilde's \cite{DeW2} popular work on the closed graph theorem led to webbed spaces being a common version of that theorem, Valdivia  \cite{QLB} used webbed space structures to prove  distinct closed graph and lifting theorems for locally convex spaces that were defined as quasi - (LB) spaces.  Valdivia's results seemingly attracted less attention than might have been expected, perhaps as `yet another closed graph theorem'.  In fact, Valdivia's  methods and  results of \cite{QLB} still lead an active life in research.  In particular, the general structure of quasi- (LB) spaces is used quite successfully in areas such as compact resolutions in topological groups \cite{Gab}, descriptive topology \cite{Ferr0},  the Fr\'{e}chet - Urysohn property in various contexts \cite{cks},  and $K$ analyticity in spaces of continuous functions, \cite{Tka}, to name a few. In this paper, we appeal to W. Robertson's \cite{WR} original definition of webs that facilitates convergence of sequences in non-metrizable locally convex spaces so as to extend the lifting theorem of \cite{QLB} from metrizable strictly barrelled   locally convex  spaces as the range, to those having sequential webs of \cite{TG}.  The collection of locally convex spaces with sequential webs includes those that are metrizable, and also includes some collections of locally convex spaces  that are not necessarily metrizable, such as any sequentially retractive locally convex inductive limit; c.f. \cite{TG}.  Thus, for example, the space $\mathcal{D}$ of test functions from the classical Schwartz distribution theory is allowable here. In addition, the sequential web condition does not require that the spaces be strictly barrelled or even barrelled. Examples of range spaces are  provided, including that of  the space  $\mathcal{D}_{\Gamma}^{\prime}$  of distributions having their wevefront sets in a specified closed cone $\Gamma$, equipped with a locally convex topology, (denoted as the normal topology, \cite[p. 1350]{dab}).  This space, and its dual space, $\mathcal{E}_{\Lambda}^{\prime}$  have become important in mathematical physics, due to their role in the formulation of quantum field theory in curved spacetime \cite[p.1345]{dab}.  We  observe that Valdivia's closed graph theorem is also applicable to domain spaces that need not be strictly barrelled, and could be applied  to locally convex spaces that need not even be barrelled. Thus,  Valdivia's results from \cite{QLB} are deeper and more applicable than at first glance.    \\

 Topological vector spaces having completing webs are commonly used for proving closed graph and open mapping theorems,   \cite{Bam},  \cite{Blon},  \cite{Casas1}, \cite{Casas2}, \cite{Casc}, \cite{DeW1}, \cite{Fern},  \cite{Gach}, \cite{Gar},  \cite{Ruess1},   \cite{WR}, \cite{V1},  \cite{QLB}, \cite{Vin}, among others.  In  \cite{QLB},  strict webs, originally defined as strict r\'{e}seau absorbants  (see \cite{DeW2}),  are equipped with an  ordered structure and used to prove  substantial generalizations of the  closed graph and lifting theorems of DeWilde \cite{DeW2} using techniques that connect  back to those of Banach's original closed graph theorem proof, \cite{Ban}.    The  approach taken in  \cite{QLB}  is to assume that the linear function and strands of a web satisfy properties for  closed graph, open mapping, and lifting theorems, rather than assume the domain and range spaces be of a particular type of locally convex space. See  \cite[9.1.44, p. 346]{PCB}  for another example of this approach.    By replacing metrizability with the presence of sequential webs,  Valdivia's lifting  theorem  \cite[6.4, p.  162]{QLB}   can be generalized to  spaces that need  be neither strictly barrelled nor metrizable.  \\
 
  \noindent   Valdivia's original Lifting theorem is stated here:  \\
 
  \noindent {\bf Valdivia's Lifting theorem.}  (\cite[6.4, p. 162]{QLB}).  {\it Suppose $E$ and $F$ are Hausdorff  locally convex spaces such that $F$ is a quasi-(LB) space  and  $E$ is metrizable and strictly barrelled.  Let $T$ be a linear map from   subspace $H$ of $F$ onto $E$ for which the graph  of  $T$ is fast closed in $F \times E$.  Suppose  $\{ A_{\alpha} : \alpha\in \mathbb{N}^{\mathbb{N}}\}$ is a quasi-(LB) representation in $F$.   Define}  
  
  $$ V_{m_{1} \cdots m_{k}} = \bigcup \{T(A_{\alpha} \cap H) : \alpha \in  \mathbb{N}^{\mathbb{N}},  \alpha = (a_{n}), \, a_{n} = m_{n}, n\in \underline{k}\}.   $$
  
  \noindent  {\it Assume  there exists a sequence $(r_{k})$ of positive integers such that $\overline{V_{r_{1} \cdots r_{k}}}$ is a neighborhood of the origin in $E$ for $k=1, 2, \cdots$.  Then the following properties hold:}
   \begin{enumerate}
     \item[(a)]  {\it If $(x_{n})$ is a sequence in $E$ which converges to the origin, there exists a sequence $(u_{n})$  in $H$, fast convergent to the origin in $F$ and such that} $Tu_{n} = x_{n}, n = 1, 2, \cdots$. 
            
 \item[(b)]    {\it If  $A$  is a precompact set   in  $E$, then there exists a subset $B$ in $H$,  fast precompact in $F$ and such that} $T(B) = A$.      
      
 \item[(c)]  {\it If $(y_{n})$ is a  Cauchy sequence in $E$, there exists  a sequence $(v_{n})$ in $H$, fast convergent in $F$ and such that}  $Tv_{n} = y_{n}, n=1, 2, \cdots$.     
\end{enumerate}

 \bigskip
 
 \noindent We will abbreviate  the  above assumption that   there exists a sequence $(r_{k})$ of positive integers such that $\overline{V_{r_{1} \cdots r_{k}}}$ is a neighborhood of the origin in $E$ for $k=1, 2, \cdots$  with  the term {\it strand-neighborhood condition} and replace the assumptions that $E$ is strictly barrelled and metrizable with the less restrictive assumption that $E$ has a  sequential web.  The  main result of this paper is:  \\

  \begin{thm}  (Lifting Theorem). Suppose $E$ and $F$ are Hausdorff  locally convex spaces such that $F$ is a quasi-(LB) space  and  $E$ has a sequential web.   Let $T$ be a linear map from a fast closed linear subspace $H$ of $F$ onto $E$ for which the graph  of  $T$ is  closed in $H\times E$.  If the images under $T$  of the quasi- (LB) representation of $F$ restricted to  $H$  satisfy the strand - neighborhood condition, then:
   \begin{enumerate}
     \item[(a)]   Every convergent sequence in $E$  is the image under $T$ of a  sequence in $H$ that is fast convergent in $F$.
       
    \item[(b)]    If in  $E$  each  precompact set is contained in the absolutely convex hull of a null sequence, then  every precompact subset of  $E$ is the image under $T$ of a  subset of $H$ that is fast  precompact in $F$.   
      
 \item[(c)]  Suppose that in addition $E$ is metrizable.  Then every Cauchy sequence in $E$ is the image under $T$ of a sequence in  $H$ that is fast convergent in $F$.  

 
\end{enumerate}
 \end{thm}

\noindent    Theorem 1  generalizes  \cite[6.4, p. 162]{QLB} to range spaces $E$ that possess sequential webs but need not be metrizable and need not be strictly barrelled, as will be seen in the examples presented in section 4.  The outline of this paper is as follows:  The next section will consist of basic definitions, including those  of webs,  r\'{e}seau, quasi (LB)  spaces, strictly barrelled spaces, sequential webs,  as well as  statements of relevant  results that will be used.  Section three consists of the proof of Theorem 1 and a couple of consequences thereof.  The penultimate  section consists of distinguishing examples and basic properties of the spaces involved.  Of note is that many  natural examples of such spaces turn out to be non-barrelled, and this implies potential applications to such spaces.  This paper concludes with a few  indications regarding potential further developments.   \\

 \section{Preliminaries.}

\subsection{General}   Throughout this paper, $E = (E, \tau)$ denotes a Hausdorff  locally convex topological vector space, or {\bf space},   over $\mathbb{K} = \mathbb{R} \mbox{ or } \mathbb{C}$.  A convex, balanced set is {\bf absolutely convex}, which we will call  a {\bf disk}.   The expression $absconv(A)$ denotes  the absolutely convex hull of a set $A$.   A system of zero neighborhoods in a space $E$ with topology $\tau$ will be denoted by $\mathcal{N}_{0, E}$.  When distinction of topologies is needed, we use the notation $\tau_{E}$ for a space $E$.  The set $\{1, 2, \cdots, k\}$ is denoted by $\underline{k}$.  The closure of a set $A$ in   $(E, \tau)$ will be denoted by $\overline{A}$, and by $\overline{A}^{\, \tau_{E}}$ whenever it is necessary to clarify the topology involved. Likewise, the interior of a set $A$ is $A^{\circ}$, and $A^{\circ \, \tau_{E}}$ to clarify.   A linear subspace $H$ of a space $F$ will be denoted by $H\leq F$.  We make typical use of Banach disks and related ideas described next, which incidentally,  are originally due to Grothendieck  \cite{Groth}:  \\

\begin{defn} \begin{enumerate}  
\item Given an absolutely convex   subset $B$ of a space $E$, the linear span of $B$ is denoted by $E_{B}$, and we equip $E_{B}$ with the linear topology given by the Minkowski functional of $B$, namely, for any $x\in E_{B}$, $\mu_{B}(x) = \inf\{t > 0 : x\in tB\}$.   If additionally the disk $B$ is  bounded, then  with this topology $(E_{B}, \mu_{B})$ is a normed space and we write  $(E_{B}, ||\cdot ||_{B})$. It is easy to see that a base of zero neighborhoods for $(E_{B},  ||\cdot ||_{B})$ is given by $\{n^{-1} B : n\in\mathbb{N}\}$. The boundedness of $B$ implies that the injection $(E_{B}, ||\cdot ||_{B})   \hookrightarrow E$ is continuous.  When this normed space is complete,  $B$ is called a {\bf Banach disk}.

\item A  sequence $(y_{n})$  in $E$  is  {\bf Mackey convergent to} $y$ if converges  to $y$ in  $(E_{B}, ||\cdot ||_{B})$, for some  closed bounded disk $B$, \cite[5.1.29, p. 158]{PCB}.  In particular, $(y_{n})$ is  {\bf fast convergent},  if the disk  $B$ is a Banach disk, \cite[6.1.20, p. 171]{PCB}.    Equivalently a sequence $(y_{n})$  in $E$  is  Mackey convergent  to $y$ \cite[5.1.1, p. 151]{PCB}   if there exists a sequence $(d_{n})$ of positive numbers such that $d_{n}\rightarrow \infty$ as $n\rightarrow \infty$, and $d_{n}  (y_{n} - y) \rightarrow 0$ in $E$.  A space $E$ satisfies the {\bf Mackey convergence condition},  if any null sequence is Mackey null \cite[5.1.29, p. 158]{PCB}.    

\item     A linear subspace $H$ of a space $F$ is   {\bf fast closed} if, given a Banach disk $B$ of $F$,  $H\cap F_{B}$  is a closed subspace of $F$, cf \cite[p. 150]{QLB} .    

\item   A subset $A$ of $E$ is {\bf fast (pre) compact} if for some Banach disk $B$, $A$ is (pre) compact in $(E_{B}, ||\cdot ||_{B})$.

\end{enumerate}
\end{defn}   

\begin{defn}    An inductive limit $E = ind_{n}E_{n}$ is {\bf sequentially retractive} \cite[8.5.32, p. 295]{PCB},  if given any null sequence  in $E$, there exists $m\in \mathbb{N}$ such that the sequence is contained in and null in $E_{m}$.     \end{defn} 

\begin{defn}  (See \cite[Def 1.1, p. 52]{Dom}).  A projective limit of a sequence of strong duals of Fr\'{e}chet-Schwartz spaces (i.e., DFS spaces) is called a {\bf PLS - space}.  \end{defn}

 \noindent  So as to distinguish the results here from other concepts of ``lifting'',  we include (c.f. \cite[VI.3.5, p. 117]{DeW1}):
 \begin{defn}  Consider a linear map $T: F\rightarrow E$, for locally convex spaces $E$ and $F$.  Then $T$ is a {\bf lifting} if either of the following hold:
 \begin{enumerate} \item Every (pre-) compact subset of $E$ is the image under $T$ of a (pre-) compact subset of $F$.  
 \item Every (fast)  (pre-) compact set in $E$ is the image under $T$ of a fast (pre-) compact subset of $F$.  
 \end{enumerate}
 \end{defn}
  
\noindent General background information on topological vector spaces  can be found in \cite{PCB} or \cite{RR}.   \\

\subsection{Webs and R\'{e}seaux webs (absorbents) and a comparison}  DeWilde`s  r\'{e}seaux absorbents  are  commonly denoted with  W. Robertson's original terminology of webs \cite{WR}. The two definitions give rise to distinct representations, and because we will use both constructions, a slightly closer look at each is relevant.  Both definitions    follow.

 \begin{defn} (\cite[p. 714]{WR})   A  {\bf (Robertson) web} $\mathcal{W}$ on a space $F$ is a countable collection $\mathcal{W}$ of balanced sets,  arranged in {\it layers}, where the first layer is given by $\{W_{n_{1}} : n_{1}\in \mathbb{N}\}$ and subsequent layers, denoted by $\{W_{n_{1} n_{2} \cdots n_{k}} : k,  n_{1}, n_{2}, \cdots, n_{k} \in \mathbb{N}\}$, that satisfy the following properties:  \\

\begin{enumerate}
 \item For each $k\in\mathbb{N}$,  $$W_{n_{1} n_{2} \cdots n_{k+1}}  + W_{n_{1} n_{2} \cdots n_{k+1}}  \subset W_{n_{1} n_{2} \cdots n_{k}};     $$

\item     $$\bigcup\{ W_{n_{1}} : n_{1}\in \mathbb{N} \} \mbox{ is absorbing in }  F, $$

\item    For each $k\in\mathbb{N}$, $$\bigcup \{W_{n_{1} n_{2} \cdots n_{k+1}} : n_{k+1} \in \mathbb{N}\} \mbox{ is absorbing in }   W_{n_{1} n_{2} \cdots n_{k}}.  $$

 \item   Given $(n_{1}, n_{2}, \cdots) = (n_{k}) \in \mathbb{N}^{\mathbb{N}}$, the sequence  $(W_{n_{1}}, W_{n_{1} n_{2}}, \cdots)$ of sets  from $\mathcal{W}$ is called a {\bf strand}, and is denoted by $(W_{k})$, that is, 
 
 $$  (W_{k})_{k\in\mathbb{N}} = (W_{k}) = (W_{n_{1} n_{2} \cdots n_{k}}).  $$
 \noindent    In particular, when all subsets of a web are absolutely convex  then for  any strand $(W_{k})$ of $\mathcal{W}$, item (1) becomes:
 
\begin{equation}\tag{*} (\forall k\in\mathbb{N}) \; \; \; \; \;  W_{k + 1} \subset \frac{1}{2} W_{k}.  \end{equation}
 \end{enumerate}
 \end{defn}
 
 \begin{defn} (\cite[p. 14]{DeW2}, \cite{DeW1})  A {\bf r\'{e}seaux (absorbents) web}  on a  space $F$ is a countable collection $\mathcal{R}$ of sets, denoted by   $\{A_{n_{1} n_{2} \cdots n_{k}} : k,  n_{1}, n_{2}, \cdots, n_{k} \in \mathbb{N}\}$,  satisfying the following properties:  
 
 \begin{enumerate}

\item     $$\bigcup\{ A_{n_{1}} : n_{1}\in \mathbb{N}\} = F, $$

\item    For each $k\in\mathbb{N}$, $$\bigcup \{A_{n_{1} n_{2} \cdots n_{k+1}} : n_{k+1} \in \mathbb{N}\} =   A_{n_{1} n_{2} \cdots n_{k}}.  $$
  \end{enumerate}
  \end{defn}
  
  \bigskip
  
  \noindent  The following is the completeness property usually assumed for webs or  r\'{e}seaux webs, on a  space $F$.
  
 \begin{defn}  \begin{enumerate}  \item  (\cite[p.715]{WR})    A web $\mathcal{W}$  in  $F$   is  of {\bf type (c)}  if, given any strand $(W_{k})$ of $\mathcal{W}$  and any sequence $(x_{k})$ with $x_{k}\in W_{k}$   for each $k\in\mathbb{N}$, the series $\sum_{k=1}^{\infty} x_{k}$  converges in $F$.  If in addition the  convergent series satisfy  $\sum_{r=k+1}^{\infty} x_{r} \in W_{k-1}$ for each $k\geq 2$, then $\mathcal{W}$ is {\bf tight}.  
 
 \item   (\cite[p. 14 - 15]{DeW2}, \cite[p. 48 - 49]{DeW1})    A r\'{e}seaux web $\mathcal{R}$   in  $F$ is of {\bf type $\mathcal{C}$}  if  given any sequence $(n_{k})  \in \mathbb{N}^{\mathbb{N}}$, and any  $x_{k}\in  A_{n_{1} n_{2} \cdots n_{k}}$ for each  $k\in\mathbb{N}$,  there exists a sequence $(\lambda_{k})$ of positive numbers, such that the series $\sum_{k=1}^{\infty} \mu_{k} x_{k}$ converges in $F$ for all $\mu_{k}$ for which $0\leq |\mu_{k}| \leq \lambda_{k}$, $k \in\mathbb{N}$.    If in addition the  convergent series satisfy   $\sum_{k=p}^{\infty} \mu_{k} x_{k} \in A_{n_{1} n_{2} \cdots n_{p}}$, for any $p\in\mathbb{N}$, then  $\mathcal{R}$ is {\bf strict}.   
 
   \end{enumerate}
 \end{defn}
  
  \noindent  It is assumed here that  all sets in a web or r\'{e}seaux web are absolutely convex.   As observed in \cite[p. 725]{WR},  a  web of type (c)  is equivalent to a r\'{e}seaux  of type $\mathcal{C}$, and a tight web is equivalent to a strict web.  A space with a completing web or  r\'{e}seaux web is typically referred to as a {\bf webbed} space, and spaces with tight  webs or strict r\'{e}seaux webs are referred to as  {\bf strictly webbed} spaces or having a {\bf strict web}.  This is how we will generally refer to webs in the sequel.  We need a few more properties of webs:  
  
 \begin{defn}   A web $\mathcal{W}$ on  a space $E$ is:   
  \begin{enumerate}
 \item {\bf  Compatible}  \cite[p. 156]{RR} if, given any zero neighborhood $V$ in $E$, and any strand $(W_{k})$ of $\mathcal{W}$, there exists $K\in \mathbb{N}$ such that $W_{K} \subset V$.  Hence, in this case, $(\forall k\geq K) \, W_{k} \subset V$. 
 
 \item {\bf Sequential} \cite[Def. 6, p. 475]{TG} if it is compatible, and given any null sequence $(x_{n})$ in $E$, there is a strand $(W_{k})$ of $\mathcal{W}$ for which, given any  $k\in \mathbb{N}$, there is $N_{k}\in\mathbb{N}$  such that  $x_{n}\in W_{k}$, for all $n\geq N_{k}$.  In this construction, one assumes the property of strands $(W_{k})$  of W. Robertson's webs c.f. \cite[p. 156]{RR} as expressed in (*) of Definition 5. 
 \end{enumerate}
  \end{defn}
  
  \noindent  If $E$ is any metrizable  space, and $\{V_{k} : k\in\mathbb{N}\}$ is a system of absolutely convex zero neighborhoods such that $(\forall k\in\mathbb{N}) \; 2V_{k+1} \subset V_{k}$, then $\mathcal{W} = \{V_{k} : k\in\mathbb{N}\}$ is a sequential web.  In particular, the $k$-th layer of  $\mathcal{W}$ consists of the set $V_{k}$.    General references for webs on topological vector spaces include \cite{RR}, \cite{DeW1}, and \cite{WR}. 
    
 \subsection{Ordered webs and  quasi-(LB) spaces.}  

\begin{defn}  \cite[p. 150]{QLB}  A web $\{W_{n_{1} n_{2} \cdots n_{k}} :  k, n_{1}, n_{2}, \cdots, n_{k} \in \mathbb{N}\}$ is {\bf ordered} if, given arbitrary natural numbers $k, r_{1}, r_{2}, \cdots, r_{k}$, and $s_{1},  s_{2}, \cdots, s_{k}$, if $r_{j} \leq s_{j}$ for all $j\in\underline{k}$, then 

$$   W_{r_{1} r_{2} \cdots r_{k}} \subset W_{s_{1} s_{2} \cdots s_{k}}.    $$
\end{defn}

 \begin{defn}  \cite[ p. 151]{QLB}  Given $\alpha  = (a_{n})\in \mathbb{N}^{\mathbb{N}}$ and $\beta  = (b_{n})\in \mathbb{N}^{\mathbb{N}}$, one puts $\alpha \leq \beta$ if $a_{n} \leq b_{n}$ for every $n\in \mathbb{N}$.  A {\bf quasi-(LB) representation} on a space $F$  is  a family of Banach disks  $\{A_{\alpha} : \alpha \in  \mathbb{N}^{\mathbb{N}}\}$ such that:

 \begin{enumerate}
 
 \item $\bigcup \{ A_{\alpha} : \alpha \in  \mathbb{N}^{\mathbb{N}}\} = F$;    
 
 \item For any $\alpha, \beta \in \mathbb{N}^{\mathbb{N}}$, if $\alpha \leq \beta$, then $A_{\alpha} \subset A_{\beta}$.
 \end{enumerate}
 For short, we call such an $F$ a {\bf quasi- (LB) space}.    For a given   $\alpha \in  \mathbb{N}^{\mathbb{N}}$, we denote by $F_{A_{\alpha}}$ the corresponding Banach space  $(F_{A_{\alpha}}, || \cdot ||_{A_{\alpha}})$.  
\end{defn}

\noindent  The connection between quasi- $(LB)$ spaces and webs that we will make use of is given by the following.

 \begin{thm}    \cite[4.1, p. 153]{QLB} A locally convex space is a quasi-(LB) space if, and only if, it has an absolutely convex, ordered,  strict web.    In that case, the strands are described by:  Given a sequence $(m_{n})\in \mathbb{N}^{\mathbb{N}}$,
 \begin{equation}  U_{k} \equiv  U_{m_{1} m_{2} \cdots m_{k}} = \bigcup \{ A_{\alpha} : \alpha = (a_{n})\in\mathbb{N}^{\mathbb{N}}, a_{n} = m_{n}, n\in\underline{k}\}.  
  \end{equation}
\end{thm}

\noindent  Observe that  $U_{m_{1}} \supset  U_{m_{1} m_{2}} \supset \cdots$, for each $(m_{n})\in \mathbb{N}^{\mathbb{N}}$.  \\
     The proofs in \cite{QLB} use ordered  (strict)  r\'{e}seaux web structures  as given in Definitions 6 and 9.    The parts of  proofs here that correspond to sequential webs are  based on the (tight) web structures of Definition 5, where sequential convergence is easy to manage.   The following simple example illustrates the difference between an ordered, tight web and and ordered strict  r\'{e}seaux web, as being equivalent, but having distinct representations on a space.  

\begin{ex}  Tight ordered webs and strict ordered  r\'{e}seaux webs on any nontrivial $(LB)$ - space.
\end{ex}
\noindent  {\it Details}.   Let  $E = indlim_{n}(E_{n})$, where, for each $n\in \mathbb{N}$, $B_{n}$ denotes the closed unit ball of the Banach space $E_{n}$.  Multiplying by appropriate scalars if necessary, we can assume that for each  $n\in\mathbb{N}$, one has $B_{n}\subset B_{n+1}$.    An ordered,  tight web on $E$ is described by setting the $k$th layer to be $\{ 2^{-k} B_{n} : n\in\mathbb{N} \}$, for each $k\in\mathbb{N}$.  On the other hand, an ordered, strict r\'{e}seaux web is described as follows:  Certainly, $E = \bigcup_{n=1}^{\infty} n B_{n}$, for the first layer.  For the second layer, one writes each $nB_{n}$ as an increasing countable union of Banach disks. Subsequent layers are obtained by repeating this process. Using values of the form $2^{j}$ for appropriate values of $j \in\mathbb{Z}$, one can start with an ordered tight web on $E$ and construct a strict  r\'{e}seaux web representation, and visa versa.   

\begin{defn}  \cite[p. 160]{QLB}   A space is {\bf strictly barrelled} if for every ordered, absolutely convex web on the space, there exists at least one strand for which the closure of each member  is a  zero neighborhood.  \end{defn}

\noindent It is easy to see that every strictly barrelled space is barrelled.  A barrelled space that is not strictly barrelled will be given in section 4.  \\

\noindent    The following definition  is intended to express the general assumption that linear images or linear  inverse images of strands of a certain  ordered strict web structure  satisfy a condition of closures being zero neighborhoods.  See  \cite[p. 155, 157]{QLB}.

  \begin{defn}   Consider spaces $E$ and $F$, and a linear map $T: E\rightarrow F$.   Assume $\{ A_{\alpha} : \alpha \in  \mathbb{N}^{\mathbb{N}}\}$ is a quasi - (LB) representation in $F$. Given $k, m_{1}, m_{2}, \cdots, m_{k} \in \mathbb{N}$, put
  
   \begin{equation} U_{m_{1} m_{2} \cdots m_{k}} = T^{-1}\left( \bigcup \{ A_{\alpha} : \alpha = (a_{n})\in\mathbb{N}^{\mathbb{N}}, a_{n} = m_{n}, n\in\underline{k}\}\right).  
  \end{equation}
  
  The triple $(E,F,T)$ satisfies a {\bf strand - neighborhood condition}  if there exist a sequence $(r_{k})\in  \mathbb{N}^{\mathbb{N}}$   such that for  $U_{k} \equiv  U_{r_{1} r_{2} \cdots r_{k}}$,  

 $$  \overline{U_{k}} =   \overline{U_{r_{1} r_{2} \cdots r_{k}}} \in \mathcal{N}_{0, E},  $$
 
\noindent  for each $k\in \mathbb{N}$. \\
 
\noindent   For  the case of of a linear $T : F \rightarrow E$  such that $T$ maps a subspace $H$ of $F$ onto $E$,  one defines  \cite[p. 157]{QLB}:  Given $n, m_{1}, m_{2}, \cdots, m_{n} \in \mathbb{N}$, put
  
   \begin{equation}  V_{m_{1} m_{2} \cdots m_{n}} = \bigcup \{ T(A_{\alpha} \cap H): \alpha = (a_{n}) \in\mathbb{N}^{\mathbb{N}}, a_{n} = m_{n}, n\in\underline{k}\}.  
  \end{equation}

\noindent We assume there exists $(r_{k})\in  \mathbb{N}^{\mathbb{N}}$   such that for  $V_{k} \equiv  V_{r_{1} r_{2} \cdots r_{k}}$,  

 $$  \overline{V_{k}} =    \overline{V_{r_{1} r_{2} \cdots r_{k}}} \in \mathcal{N}_{0, E},  $$
  \noindent  for each $k\in \mathbb{N}$. \\
   \end{defn}
   
 \noindent    We state   Valdivia's closed graph theorem  below, as its use will be important in the proof of the main result. We use the notation of Definition 12.    \\
   
\noindent {\bf Valdivia's Closed Graph Theorem}  (\cite[5.11, p. 157, 6.1.5, p. 163]{QLB}).  A linear map from a Hausdorff locally convex space to a quasi -(LB) space having a closed graph is continuous if the strand neighborhood condition is satisfied. In particular,  we have, in $E$, 
 \begin{enumerate}
 
 \item $$ \emptyset \neq \overline{U_{k}}^{\circ \, \tau_{E}} =  \overline{U}^{\circ \, \tau_{E}}_{r_{1} r_{2} \cdots r_{k}}  \subset U_{r_{1} r_{2} \cdots r_{k}},  \; k\in\mathbb{N},  $$   
 \item  $T$ is continuous.
 \end{enumerate}
  
 \noindent  One more result is needed in the proof of Theorem 1:
 
 \begin{lem}  \cite[Lemma 3, p. 60]{Rest}  Let $T:F\rightarrow E$ be linear, for topological vector spaces $E$ and $F$.  Consider   
 
 $$  F \stackrel{\varphi}{\rightarrow} F/T^{-1}(\{0\}) \stackrel{\tilde{T}}{\rightarrow}  E, $$
 
 \noindent where $\varphi$ is the canonical linear quotient map, and $\tilde{T}$ is the linear map such that $T = \tilde{T}\circ \varphi$.  If $G_{T}$ denotes the graph of $T$, then $G_{T}$ is closed if and only if, $G_{\tilde{T}}$ is closed.  \\
 \end{lem}

  
\section{Proof of the main result.}

\noindent {\it Proof of Theorem 1}.  Our methods here are inspired by those in \cite{Casas2}, \cite{Casc},  \cite{Fern}, and  \cite{QLB}.  For the sake of completeness, full details are given here.   It suffices to prove  (a) and (b), because the proof of part (c) only requires the observation that any metrizable space has a sequential web and is otherwise  exactly the same as in \cite[5.6, p. 158]{QLB}.  \\

(a).  The linear map $T$ is not assumed to be continuous, and the closed graph theorem does not apply to $T:F \rightarrow E$ (the assumptions are reversed).  Nevertheless, the arrival at the conclusion  can be obtained by taking a detour through the quotient  space $H/T^{-1}(\{0\})$ in such a way as to apply Valdivia's closed graph theorem. This is done in Step 2.  After that, we apply the assumption of  a sequential web to  prove the existence of the desired sequence. By way of the linearity of topologies and functions, it suffices to consider null sequences in our work.   \\

 \noindent  {\it Step 1}.  {\it Setting up  spaces and linear maps}.  Let $T$ be a linear map from a fast closed linear subspace $H$ of $F$ onto $E$, for which the graph $G_{T}$ of $T$ is  closed in $H\times E$.  Let   $\{ A_{\alpha} : \alpha \in \mathbb{N}^{\mathbb{N}} \}$  be a  quasi-(LB) representation on  $F$.  Suppose $E$ has a sequential web denoted by  $\mathcal{Z}$, for which we denote strands by $(Z_{k})_{k\in\mathbb{N}} = (Z_{k})$.  Consider any  null sequence $(y_{n})$  in $E$.   We wish to show there exists  $\beta \in \mathbb{N}^{\mathbb{N}}$  and a sequence $(u_{n})$ in $F_{A_{\beta}} \cap H$ such that $(u_{n})$ is fast convergent to zero in $F$, and $T(u_{n}) = y_{n}$, for all $n\in \mathbb{N}$.     Using previous notation, we have   
 
\begin{equation}    V_{k} \equiv   V_{m_{1} \cdots m_{k}} = \bigcup \{T(A_{\alpha} \cap H) : \alpha \in  \mathbb{N}^{\mathbb{N}},  \alpha = (a_{n}), \, a_{n} = m_{n}, n\in \underline{k}\}.    \end{equation} 
 
\noindent  By the strand-neighborhood assumption,   there exists a sequence $(r_{k}) \in \mathbb{N}^{\mathbb{N}}$ such that for all $k\in\mathbb{N}$,   $\overline{V}_{r_{1} \cdots r_{k}} \in\mathcal{N}_{0, E}$.  \\


  We proceed with a  factoring through the null space $T^{-1}(\{0\})$, as follows:  Denote by $\varphi$ the canonical quotient map from $H$ onto $H/T^{-1}(\{0\})$.  Next, let $S : E {\longrightarrow} H/ T^{-1}(\{0\})$ be the bijective  linear map such that $T = S^{-1} \circ \varphi$.   Observe that for all $\alpha \in  \mathbb{N}^{\mathbb{N}}, \, \varphi(A_{\alpha} \cap H) = S(T(A_{\alpha} \cap H))$.    Because $H$ is a fast closed linear subspace of $F$, for every $\alpha \in  \mathbb{N}^{\mathbb{N}}$, $ \overline{H\cap F_{a_{\alpha}}}^{F} = H\cap F_{a_{\alpha}}$, hence, $\{ A_{\alpha} \cap H : \alpha \in \mathbb{N}^{\mathbb{N}}\}$ is a quasi- (LB) representation in  $H$.  By the continuity of $\varphi$ from $H$ onto $H/T^{-1}(\{0\})$,  an application of  \cite[3.1, p. 151]{QLB}  guarantees that $\{ \varphi(A_{\alpha} \cap H) : \alpha \in \mathbb{N}^{\mathbb{N}}\}$  is a quasi- (LB) representation in  $H/T^{-1}(\{0\})$.  

\[
\begin{tikzcd}[column sep=huge,row sep=huge]
H\leq F \arrow[r,"T = \, S^{-1} \circ \; \varphi"] \arrow[dr,swap,"\varphi"] &
  E \arrow[d,shift left=.75ex,"S"] \\
& H/ T^{-1}(\{0\})  \arrow[u,shift left=.75ex,"S^{-1}"]
\end{tikzcd}
\] 
 \bigskip
 
\noindent {\it Step 2}.  {\it Applying Valdivia's closed graph theorem to the linear map $S : E \rightarrow H/T^{-1}(\{0\})$}.  We  apply Lemma 1 to $T$,  using $H$ in place of $F$, and $\tilde{T} = S^{-1}$:   The assumption that the graph of $T$ is closed in $H\times E$ implies  that the graph  of $S^{-1}$ is closed in $H/T^{-1}(\{0\}) \times E$.  Because $S$ is bijective, we conclude that the graph of $S$ is closed as well.   Consider the linear map $S : E \rightarrow H/T^{-1}(\{0\})$ with  $H/T^{-1}(\{0\})$ being  a quasi- (LB) space with  representation $\{ \varphi(A_{\alpha} \cap H) : \alpha\in\mathbb{N}^{\mathbb{N}}\}$.   If we   verify the strand - neighborhood condition  $S$, then Valdivia's closed graph theorem can be applied to $S$.  Given

\begin{flalign} C_{m_{1} \cdots m_{k}} &  = \bigcup \{\varphi(A_{\alpha} \cap H) : \alpha \in  \mathbb{N}^{\mathbb{N}},  \alpha = (a_{n}), \, a_{n} = m_{n}, n\in \underline{k}\} \nonumber  \\
                                                                & =  \bigcup \{S(T(A_{\alpha} \cap H) ): \alpha \in  \mathbb{N}^{\mathbb{N}},  \alpha = (a_{n}), \, a_{n} = m_{n}, n\in \underline{k}\},  \nonumber   
    	\end{flalign}
	
\noindent  put  $$U_{m_{1} \cdots m_{k}} =S^{-1}(C_{m_{1} \cdots m_{k}}).  $$

\noindent  By (3) of Definition 12, we obtain 

\begin{flalign}  S^{-1}(C_{m_{1} \cdots m_{k}}) = S^{-1}\left( \bigcup \{\varphi(A_{\alpha} \cap H) : \alpha \in  \mathbb{N}^{\mathbb{N}},  \alpha  = (a_{n}), \, a_{n} = m_{n}, n\in \underline{k}\} \right)   \nonumber \\
 = \bigcup \{T(\varphi(A_{\alpha} \cap H) : \alpha \in  \mathbb{N}^{\mathbb{N}},  \alpha = (a_{n}), \, a_{n} = m_{n}, n\in \underline{k}\}  = V_{k}.  \nonumber 	\end{flalign}

\noindent  Therefore, with the  same $(r_{k}) \in  \mathbb{N}^{\mathbb{N}}$ as in the hypotheses,    

 $$  \overline{V}^{\, \tau_{E}}_{k} = \overline{V}^{\, \tau_{E}}_{r_{1}  r_{2} \cdots r_{k}} \in \mathcal{N}_{0, E}, $$

\noindent for all $k\in\mathbb{N}$.  It now follows that the strand-neighborhood condition is satisfied for $S$, and we invoke  Valdivia's closed graph theorem to conclude: 

 \begin{enumerate}
 \item[(i).] $ \emptyset \neq \overline{V_{k}}^{\circ \, \tau_{E}} =  \overline{V}^{\circ \, \tau_{E}}_{r_{1} r_{2} \cdots r_{k}}  \subset V_{r_{1} r_{2} \cdots r_{k}},  \; k\in\mathbb{N}$,     \\
 
 \item[(ii).]  $S$ is continuous.
 \end{enumerate}

\noindent  In particular, each $V_{r_{1} r_{2} \cdots r_{k}}$ is a zero neighborhood in $E$.  \\



 \noindent   {\it Step 3}.   It is now time to apply the assumption of the sequential web for the  null sequence $(y_{n})$  in $E$.  Let  $(Z_{k})$ of $\mathcal{Z}$ be as in  the hypotheses.  Given $k\in\mathbb{N}$, find $j_{k} > k$ such that $2^{-j_{k}} < k^{-1}$, and observe that
 
 $$ Z_{j_{k} + k} \subset 2^{-1} Z_{j_{k} + k-1} \subset  2^{-2} Z_{j_{k} + k -2} \subset \cdots \subset 2^{-k} Z_{k} \subset k^{-1} Z_{k},   $$
 
 \noindent  as a consequence of  the property $Z_{k+1}\subset \frac{1}{2}  Z_{k}$.  Next, the sequential property of $\mathcal{Z}$ implies the existence of  a positive integer $N_{j_{k}}$ such that for all $n\geq N_{j_{k}}$, $y_{n} \in Z_{j_{k} + k} \subset k^{-1} Z_{k}$.  \\

\noindent Define, for all $k  \in\mathbb{N}$, $d_{n} = k, N_{j_{k}} \leq n < N_{j_{k+1}}$.  Clearly, $\lim_{n\rightarrow \infty} d_{n} = \infty$.  We apply the compatibility of $\mathcal{Z}$ to the zero neighborhoods $V_{l} = V_{r_{1} r_{2} \cdots r_{l}}, \, l \in \mathbb{N}$, as follows.  For a fixed $l\in \mathbb{N}$, there exists $K_{l}\in\mathbb{N}$ such that for all $k\geq K_{l}$, $Z_{k}\subset V_{l}$.  Next, by the construction at the beginning of this step, there exists $N_{j_{K_{l}}} \in\mathbb{N}$ such that for all $n\geq N_{j_{K_{l}}}$, $d_{n}y_{n} \in V_{l}$.   The surjectivity of $T$ allows us to find a sequence $(w_{n})$ in $H$, such that for all $n\geq N_{1} \; w_{n} \in  H$ such that $T(w_{n}) = d_{n} y_{n}, \, N_{j_{k}} \leq n  < N_{j_{k+1}}$, for each $k\in\mathbb{N}$.  By   \cite[5.21 (i) p.43]{Casc} with $p = 1$, there exists $\beta  \in\mathbb{N}^{\mathbb{N}}$ such that for each $n\in\mathbb{N}$,  $w_{n}\in F_{A_{\beta}}$, where $F_{A_{\beta}}$ is the Banach space corresponding to the Banach disk $A_{\beta}$, having $\{ m^{-1} : m\in\mathbb{N}\}$ as a base for a system of zero neighborhoods. Put  $u_{n} = \frac{1}{d_{n}} w_{n}$.  Then  for all $n\geq N_{j_{K_{1}}}, w_{n} \in A_{\beta}$, hence, $u_{n} \in  \frac{1}{d_{n}} A_{\beta}$,    and because $\lim_{n\rightarrow \infty} d_{n} = \infty$,    $(u_{n})$ is fast convergent to zero in  $F_{A_{\beta}}$.  Finally,  $u_{n}\in H\cap  F_{A_{\beta}}$, and 
  
  $$  T(u_{n}) = S^{-1}\left(\varphi(\frac{1}{d_{n}} w_{n})\right) = \frac{1}{d_{n}} d_{n} y_{n} = y_{n},   $$
  
  \noindent  finishing the proof  of part (a). \\
    

   (b).   Let  $E$ be  any space such that each  precompact set is contained in the absolutely convex hull of a null sequence.  The goal is to  prove that for any precompact  $B\subset E$ , there exists $\beta \in \mathbb{N}^{\mathbb{N}}$ and a set $M\subset F_{A_{\beta}} \cap H$ such that $M$ is precompact in $F_{A_{\beta}}$, and $T(M) = B$.    By the assumption on precompact sets, there is a null sequence $(y_{n})$ in $E$ such that $\overline{absconv\{y_{n}\}}^{\tau_{E}} \supset 2B$.    Let $\beta \in\mathbb{N}^{\mathbb{N}}$ and   $(u_{n})$ in $H\cap F_{A_{\beta}}$ corresponding to $(y_{n})$, be  as in the proof of part (a).  Define $D =  \overline{absconv\{\frac{1}{2} u_{n} : n\in  \mathbb{N} \}}^{||\cdot ||_{A_{\beta}} }$.   Because the subspace $H$ is fast closed in $F$, and the injection $(F_{A_{\beta}}, ||\cdot ||_{A_{\beta}})   \hookrightarrow F$ is continuous, we obtain:
   
   $$ D  \subset  \overline{H\cap F_{A_{\beta}}}^{||\cdot ||_{A_{\beta}} } \subset \overline{H\cap F_{A_{\beta}}}^{F} = H\cap F_{A_{\beta}} \subset H.  $$
      
 \noindent  In particular, intersecting with $H$, the following is continuous:
 
  \begin{equation}\tag{**}  H\cap F_{A_{\beta}} \hookrightarrow H  \stackrel{\varphi}{\rightarrow}  H/T^{-1}(\{0\}). \end{equation}

 \noindent    That  $\left(F_{A_{\beta}}, || \cdot||_{A_{\beta}}\right)$ is a Banach space implies  $D  \subset F_{A_{\beta}}$ is compact.  The continuity of the linear mappings in (**) ensures that $\varphi(D)$ is compact in $H/T^{-1}(\{0\})$.  Next, the graph of $T$ being closed in $H\times E$ tells us that $H/T^{-1}(\{0\})$ is Hausdorff (\cite[prop. V.1, p.77]{RR}).  Thus,  $\varphi(D)$ is closed in $H/T^{-1}(\{0\})$.  We use the continuity of $S$ from part (a) to observe that the set $S^{-1}(\varphi(D))$ is closed in $E$. Meanwhile, $S^{-1}(\varphi(D)) = T(D) \supset B$, which means $T^{-1}(B) \subset D$. In order to ensure we capture exactly $B$, we intersect; that is, we  let $M =  T^{-1}(B) \cap D$.    This finishes the proof.  \/ \/ $\Box$  \\

   
\noindent  The assumption on $E$ in part (b) includes metrizable spaces by way of a variant of the Banach - Dieudonn\'{e} Theorem (\cite[10.(3)p. 273]{Ko}).  As it also follows, this collection of spaces includes any sequentially retractive  inductive limit of metrizable spaces, which of course, need not be metrizable:  Witness any strict (\cite[0.3.1, p. 2]{PCB}) proper $(LB)$ space \cite[8.5.18, p. 288]{PCB}, or the space $\mathcal{D}$ of test functions from distribution theory.    Moreover, we will see (Example 3 in the next section) that such  spaces are not strictly barrelled. At the same time, part (c) does not assume the space $E$ is strictly barrelled.   The following corollary   indicates that Theorem 1 here properly extends Valdivia's lifting theorem to more general spaces, including those often seen in applications.  \\


\begin{cor} Let $E = \lim_{n}(E_{n})$ be a sequentially retractive inductive limit of metrizable spaces, with quasi-(LB) space $F$, $H$, and $T$ as in  Theorem 1.   If the images of the quasi- (LB) representation of $H$  under $T$ satisfy the strand - neighborhood condition, then  every convergent sequence in $E$  is the image under $T$ of a  sequence in $H$ that is fast convergent in $F$, and  each  precompact subset of  $E$ is the image under $T$ of a  subset of $H$ that is fast  precompact in $F$.   \end{cor} 

\noindent   {\it Proof}: It suffices to prove that such a space $E$ has a  sequential web.  That proof is found in  \cite[Prop. 9, p. 476]{TG}.     \/ \/ \/ $\Box$.   \\

 \noindent  For the corollary below, we consider the space $\mathcal{D}_{\Gamma}^{\prime}$  of distributions having their wavefront sets in a specified closed cone $\Gamma$, equipped with its  normal topology, as defined in \cite[p. 1350]{dab}.   Apart from its importance in distribution theory, the space  $\mathcal{D}_{\Gamma}^{\prime}$ is interesting in its own right, due to the combination of topological vector space properties it does, and does not, possess.  A lucid description of several such properties can be found in \cite{dab}.  A few such properties worth mentioning here are that  $\mathcal{D}_{\Gamma}^{\prime}$  is: a normal space of distributions,  complete, nuclear, but is neither metrizable nor bornological (its strong dual, $\mathcal{E}_{\Lambda}^{\prime}$, is an incomplete inductive limit). Moreover, $\mathcal{D}_{\Gamma}^{\prime}$ is not barrelled  \cite[Cor. 36, p. 1376]{dab}.    Nonetheless, it turns out that $\mathcal{D}_{\Gamma}^{\prime}$ satisfies the properties of the range space of Theorem 1:   \\

\begin{cor}  Suppose  $F$ is a quasi-(LB) space  and  $T$ is a linear map from a fast closed subspace $H$ of $F$ onto $\mathcal{D}_{\Gamma}^{\prime}$, for which the graph  of  $T$ is  closed in $H\times \mathcal{D}_{\Gamma}^{\prime}$.  If the images of the quasi- (LB) representation of $H$  under $T$ satisfy the strand - neighborhood condition, then for any sequence $(y_{n})$ converging to 0 in $\mathcal{D}_{\Gamma}^{\prime}$, there exists  $\beta \in \mathbb{N}^{\mathbb{N}}$  and a sequence $(u_{n})$ in $F_{A_{\beta}} \cap H$ such that $(u_{n})$ is fast convergent, and $T(u_{n}) = y_{n}$, for all $n\in \mathbb{N}$, where $A_{\beta}$ is from the quasi - (LB) space representation on $F$. 
\end{cor}  

\noindent   {\it Proof}:  We only need to prove that $\mathcal{D}_{\Gamma}^{\prime}$ has a sequential web.  As a projective limit of a sequence of strictly webbed spaces, $\mathcal{D}_{\Gamma}^{\prime}$ is strictly webbed by  hereditary properties, namely, that the strong dual of a Fr\'{e}chet space as well as  projective limits of sequences of such spaces are strictly webbed, \cite{RR}.  Moreover, $\mathcal{D}_{\Gamma}^{\prime}$ satisfies the Mackey convergence condition, \cite[Lemma 21, p. 1363]{dab}, which allows us to apply  \cite[Thm. 18, p. 481]{TG}:  a strictly webbed (locally) complete space satisfying the Mackey convergence condition has a sequential web.    \/ \/ \/ $\Box$.   \\

  \section{Examples} 
  
 \noindent   In   \cite[Sect. 6, p. 160]{QLB},  the domain space of the closed graph theorem   is assumed to be strictly barrelled. Likewise, in  the lifting theorem  \cite[5.12, p. 157, 6.4, p.162]{QLB}   the range space is assumed to be  metrizable and strictly barrelled.   The following examples reveal that the collection of domain spaces in Valdivia's closed graph theorem (respectively,  range spaces of Theorem 1), are in fact strictly larger than the collection of  strictly barrelled (respectively, metrizable strictly barrelled) spaces.  \\
 
 \begin{ex}  Some spaces having a sequential web that are not strictly barrelled: 

\begin{enumerate}

\item[(a)]  Any metrizable  space that is not strictly barrelled.  
\item[(b)]  The space $(l_{1}, \sigma(l_{1}, l_{\infty}))$.  

\end{enumerate}
 \end{ex}

\noindent {\it  Details}:  (a):  Obvious. \\

 (b):  Thanks to Shur's theorem (c.f. \cite[p.85]{Diest}), weakly convergent sequences are norm convergent in $l_{1}$.  Thus, a sequential web is given by   $\mathcal{Z} = \{2^{-n} B: n \in \mathbb{N}\}$, where $B$ is the closed unit ball of $l_{1}$. The space $(l_{1}, \sigma(l_{1}, l_{\infty}))$ is, of course, neither metrizable nor barrelled.  Incidentally,  for any $1<p, q <\infty$ as conjugate exponents, $(l_{p}, \sigma(l_{p}, l_{q}))$ is neither (strictly)  barrelled, nor does it have a sequential web.    \\

\noindent  The next two examples show that the properties of having a sequential web and being strictly barrelled are distinct.  In addition, there are barrelled spaces that have sequential webs but are not strictly barrelled.  \\

\begin{ex}  Barrelled spaces having a sequential web that are  not strictly barrelled. \end{ex}  

\noindent {\it  Details}:    Let $E$ be any proper sequentially retractive inductive limit of metrizable  barrelled spaces $E_{n}$.  As an inductive limit of barrelled spaces, $E$ is also barrelled.   Consider the web formed by the strands consisting of sets of the form $\{2^{-k}V_{k,n} : k\in \mathbb{N}\}$, where, for each fixed $n\in\mathbb{N}$,  $\{V_{k,n}: k\in\mathbb{N}\}$ is a decreasing sequence of absolutely convex sets forming a system  of zero neighborhoods  in $E_{n}$ such that for every $k\in\mathbb{N}, V_{k+1,n}\subset 2^{-1}V_{k,n}$.  Sequential retractivity tells us that this construction forms  a sequential web on $E$.  This web  is also  absolutely convex,  ordered, and has the property that  the closure of every set  in every strand has  empty interior in $E$.  In particular, any strict $(LB)$- space is  barrelled, is  quasi- (LB), has a sequential web, and yet is not strictly barrelled.   \\

\noindent  For the next example, note that by  \cite[Thm 12, p. 477]{TG}, any space with a sequential web satisfies the Mackey convergence condition.   
\begin{ex}    A strictly barrelled space without a sequential web.    
\end{ex}

\noindent {\it  Details}:    Let $I$ be the index set consisting of all increasing, unbounded sequences of positive real numbers.   Put $E = \Pi_{\alpha\in I} \mathbb{R}$.  As an uncountable product of the Banach space $\mathbb{R}$, $E$ is a Baire space, by \cite[Chap. III, Sect. 7]{Bour}.    Thus, $E$ is  totally barrelled, hence strictly barrelled, by \cite[6.17, p. 160]{QLB}.  On the other hand, in \cite[ 5.1.32, p. 159]{PCB}, it is shown that $E$ does not satisfy the Mackey convergence condition, guaranteeing that $E$ cannot have a sequential web.     \\

  \noindent  We finish this section with two observations of  hereditary properties of spaces having  sequential  webs, that are relevant to Theorem 1.  In particular, in part (b) of the proposition below strict barrelledness is replaced by the property of having a sequential web, which as we have just seen, is a distinct collection of spaces. Compare to the strict barrelledness result of \cite[6.16, p. 160]{QLB}.   \\

\begin{prop}  The property of having a sequential web is preserved under:

\begin{enumerate}
\item[(a)]  Linear subspaces.
\item[(b)]   Finitely many products.
\end{enumerate} 
\end{prop}

\noindent   {\it Proof}:  For (a),   let  $\mathcal{Z}$ be a sequential web on $E$.  If $L$ is a  subspace of $E$, define $\mathcal{Z}_{L}$ to be the set of all $Z\cap L$, as $Z$ runs through the sets of $\mathcal{Z}$.  The property of having a sequential web follows immediately.  \\
Regarding (b), it suffices to prove the case for the product of two spaces $E$ and $F$, the general finite product case following by induction. Denote by  $\mathcal{Z}_{E}$ and    $\mathcal{Z}_{F}$   sequential  webs in $E$ and $F$, respectively.  For  each $k\in\mathbb{N}$, define the $k$th layer in $E\times F$ by the products of the sets in the $k$th layer in $\mathcal{Z}_{E}$ with the sets in the $k$th layer of $\mathcal{Z}_{F}$.  It is easy to verify that this web is compatible with the topology of the product $E\times F$.   Now suppose $(w_{n}) = (x_{n}, y_{n})$ is any null sequence in $E\times F$.  By assumption, there is a strand $(Z_{k})$ in  $\mathcal{Z}_{E}$ for which, given any $k\in\mathbb{N}$, there exists $N_{k}$ in $\mathbb{N}$ such that $x_{n}\in Z_{k}$ for all $n\geq N_{k}$.  Similarly,  there is a strand $(S_{k})$ in  $\mathcal{Z}_{F}$ satisfying that, given any $k\in\mathbb{N}$, there exists $M_{k}$ in $\mathbb{N}$ such that $y_{n}\in S_{k}$ for all $n\geq M_{k}$.  Putting $J_{k} = \max\{N_{k}, M_{k}\}$, for each $k\in\mathbb{N}$, one has for all $k\in \mathbb{N}$, and for all $n\geq J_{k}$: 

$$  w_{n} = (x_{n}, y_{n}) \in  Z_{k} \times S_{k}. \; \; \; \; \; \Box  $$

\section{Concluding remarks}

\noindent  Some implications for further research regarding the application of Valdivia's closed graph theorem and the lifting theorem proved here are outlined.  \\

\noindent  {\it Applying Valdivia's closed graph theorem to  spaces that are not barrelled}.  The relevant examples presented here, including non-barrelled spaces that need not have sequential webs, could appear as domain spaces for Valdivia's closed graph theorem.  \\

\noindent  {\it Applying Theorem 1 to normed or metrizable spaces that are not barrelled}.  It seems some of the easiest examples of spaces having a sequential web and are not strictly barrelled turn out to be metrizable and non- barrelled.   Far from being ``exotic'', normed and metrizable spaces that are not barrelled are rather common.  Examples of such  spaces include:  

\begin{enumerate}  \item On every infinite dimensional Banach space there exists a finer, non-barrelled norm  \cite[4.6.7 (iv), p. 131]{PCB}.  

\item  Metrizable, non-normed non-barrelled spaces can  be formed by  countable products of non-barrelled normed spaces, \cite[4.2.5, p. 103]{PCB}.

 \item From \cite[Chapt. 6]{Ferr1}, we have:  Let $l_{0}^{\infty}(\Sigma)$ denote the space of simple functions on an algebra $\Sigma$ of sets, endowed with a natural norm arising from the Minkowski functional.  Then, every separable infinite-dimensional subspace of $l_{0}^{\infty}(\Sigma)$ is a normed, non-barrelled space \cite[6.3.2, p. 127]{Ferr1}.   Furthermore, various algebras $\Sigma$ are constructed in \cite[p. 131 - 135]{Ferr1} for which $l_{0}^{\infty}(\Sigma)$ is not barrelled.   \\
\end{enumerate}

\noindent  {\it General non-barrelled spaces}.  Our results here, as well as Valdivia's closed graph theorem appear pertinent to the study of general, not necessarily metrizable, non-barrelled spaces like:   Any infinite dimensional Banach space with its weak topology, and as we have seen, as well as $\mathcal{D}_{\Gamma}^{\prime}$.    Other such spaces, such as tensor products of  spaces related to $\mathcal{D}_{\Gamma}^{\prime}$, as in  \cite{Brou}, could also be worth studying.   More examples of non-barrelled spaces and results related to  such spaces can be found in   \cite{Buch}, \cite{Ferr2},    \cite{Tsir3}, \cite{Tsir2}, \cite{Tsir1},  \cite{V2},  to mention  a few instances.  \\

\noindent  {\it Developing results in other  contexts}.  The webbed and strictly webbed space results initially proved for locally convex spaces over the fields of $\mathbb{R}$ or $\mathbb{C}$ have since been developed for general topological vector spaces, for locally convex spaces over non-Archimedean fields, for bornological vector spaces, and  for locally convex topological $\widetilde{\mathbb{C}}$ -  modules, resulting in useful applications of each those versions.  See   \cite{Bam}, \cite{Gach}, \cite{Gar}, and \cite{Vin}, respectively.    Valdivia's closed graph, open mapping, and lifting theorems have been developed for general topological vector spaces, as in  \cite{Casas2}, \cite{Casc}, \cite{Fern}.  Beyond that, no further developments have been made.  Such research  regarding Valdivia's closed graph theorem and the lifting theorem obtained here, in contexts of locally convex spaces over non-Archimedean fields,  bornological vector spaces, and  locally convex topological $\widetilde{\mathbb{C}}$ -  modules, has potential for correspondingly useful applications.

\begin{ack} I am grateful to several anonymous referees  who provided important comments that significantly  improved the quality of this paper.
\end{ack}

\end{document}